\documentclass[titlepage,
a4paper]{amsart}
\usepackage[T1]{fontenc}
\usepackage[latin1]{inputenc}
\usepackage{amsmath}
\usepackage{amssymb}

\usepackage{amsthm}

\newtheorem{definition}{Definition}

\newtheorem{rem}{Remark}
\newtheorem{theorem}{Theorem}
\newtheorem{lemma}{Lemma}

\newtheorem{corollary}{Corollary}

\newcommand{\EE}{\mathbb{E}}
\newcommand{\WW}{\mathbb{W}}

\newcommand{\RR}{\mathbb{R}}
\newcommand{\QQ}{\mathbb{Q}}
\newcommand{\NN}{\mathbb{N}}
\newcommand{\ZZ}{\mathbb{Z}}
\newcommand{\TT}{\mathbb{T}}
\newcommand{\LL}{\mathbb{L}}

\newcommand{\ualpha}{{\underline{\alpha}}}

\newcommand{\cL}{{\mathcal L}}

\newcommand{\cF}{{\mathcal F}}
\newcommand{\cM}{{\mathcal M}}

\newcommand{\cG}{{\mathcal G}}

\newcommand{\eps}{\varepsilon}
\newcommand{\aint}{{{{^*} \! \int}}}
\newcommand{\aell}{{{{^*}\ell}}}
\newcommand{\aR}{{{^*}\mathbb{R}}}
\newcommand{\aZ}{{{^*}\mathbb{Z}}}
\newcommand{\aN}{{{^*}\mathbb{N}}}
\newcommand{\aQ}{{{^*}\mathbb{Q}}}

\begin{document}

\sloppy
\pagestyle{plain}

\renewcommand{\baselinestretch}{1.5}

\title{The moment problem on the Wiener space}
\author{Frederik Herzberg}
\thanks{Abteilung f\"ur Stochastik, Institut f\"ur Angewandte Mathematik der Universit\"at Bonn, D-53115 Bonn, Germany;
Sonderforschungsbereich 611 (DFG)}
\thanks{email: herzberg@wiener.iam.uni-bonn.de}
\thanks{The author would like to thank both Professor Sergio Albeverio and Dr. Thorsten W\"ormann for helpful discussions.}

\date{}

\maketitle
\begin{abstract} 
Consider an $L^1$-continuous functional $\ell$ on the vector space of polynomials of Brownian motion at given times, suppose $\ell $ commutes with the quadratic variation in a natural sense, and consider a finite set of polynomials of Brownian motion at rational times, $f_1\left(\vec b\right),\dots,f_m\left(\vec b\right)$, mapping the Wiener space to $\RR$.

In the spirit of Schm\"udgen's solution to the finite-dimensional moment problem, we give sufficient conditions under which $\ell$ can be written in the form $\int \cdot d\mu$ for some finite measure $\mu$ on the Wiener space such that $\mu$-almost surely, all the random variables $f_1\left(\vec b\right),\dots,f_m\left(\vec b\right)$ are nonnegative.


\end{abstract}

\newpage
\renewcommand{\baselinestretch}{1.5}

\noindent

\section{Introduction} 

Consider polynomials $f_1\left(X_1,\dots,X_N\right),\dots, f_m\left( X_1,\dots,X_N\right)\in \RR\left[X_1,\dots,X_N\right]$ and a linear functional $\ell:\RR\left[X_1,\dots,X_N\right]\rightarrow \RR$, normalised in the sense that $\ell (1)=1$. When does $\ell$ represent the moments of a Borel probability distribution $\mu$ on $\RR^N$ (i.e. $\ell\left(p\left(X_1,\dots,X_N\right)\right)= \int_{\RR^N}p(y)\mu(dy)$ for all $p\left(X_1,\dots,X_N\right)\in \RR\left[X_1,\dots,X_N\right]$) such that $\mu$-almost surely, $f_1\wedge\dots\wedge f_m\geq 0$?

This question is known as the moment problem, and it was shown by Schm\"udgen that, under the condition of $\bigcap_{i=1}^m\left\{f_i\geq 0\right\}$ being compact, a sufficient and necessary criterion in the sense of this question is that $$ \forall \left(j_1,\dots,j_m\right)\in\{0,1\}^m\forall g\in \RR\left[X_1,\dots,X_N\right]\quad \ell\left(g\left(\vec X\right)^2\prod_{i=1}^m f_i\left(\vec X\right)^{j_i}\right)\geq 0$$ (using the abbreviation $\vec X:= \left(X_1,\dots,X_N\right)$). The proof has two ingredients, the first one being Haviland's Theorem (see Haviland \cite{momentH36} and Choquet \cite{momentCh69}), the second one being a result from semialgebraic geometry known as Schm\"udgen-W\"ormann Theorem (often also referred to as Schm\"udgen's Positivstellensatz --- see, for instance, Marshall \cite{momentM00}). This theorem is remarkable in its own right, and there are two fundamentally different proofs: the original functional analytic one given by Schm\"udgen \cite{momentSchm91}, and a proof based on techniques from semialgebraic geometry, most notably pre-orderings, which was found by Thorsten W\"ormann in his thesis \cite{momentW98} (see also Berr and W\"ormann \cite{momentBW01}). These results have been extended, to the case of non-compact semi-algebraic varieties by Kuhlmann and Marshall \cite{momentKM} and subsequently by Kuhlmann, Marshall and Schwartz \cite{momentKMSch}. The moment problem for signed measures has been studied, e.g., by Kounchev and Render \cite{momentKR05}. 

Conceiving of $\RR^N$ as a path space, the question of the moment problem can also be reformulated in terms of stochastic processes: Consider a stochastic process $(x_i)_{1\leq i\leq N}$ on a filtered measurable space $\left(\Omega,\left(\cF_i\right)_{1\leq i\leq N}\right)$ and suppose $\ell$ is a linear map from the vector space of polynomials in the random variables $x_1,\dots,x_N$ to $\RR$, normalised in the sense that $\ell(1)=1$. When is $\ell$ derived from a realisation of the stochastic process $(x_i)_{1\leq i\leq N}$ that is, when is there a probability measure $P$ on $\left(\Omega,\cF_N\right)$ such that $\ell\left(p\left(\vec x\right)\right)$ for all polynomials $p\in \RR\left[X_1,\dots,X_N\right]$, and the event $\bigcap_{i=1}^m\left\{f_i\left(\vec x\right)\geq 0\right\}$ has $P$-probability $1$? 

In this short paper, we shall study generalisations of the previously formulated question to other path spaces, with a particular emphasis on the Wiener space.

In order to state the main result of this paper, consider a Wiener $L^1$-continuous functional $\ell$ from $\langle\Pi_\QQ\rangle$ -- which stands for the space of all polynomials of standard Brownian motion (denoted by $\left(b_t\right)_{t\in[0,1]\cap\RR}$) and rational times, cf Definition \ref{Pi_} -- to $\RR$. Let us write $\langle\Pi_\RR\rangle$ for the space of all polynomials of Brownian motion at given times (see again Definition \ref{Pi_}).

\begin{theorem}[Corollary \ref{externalWiener} and Lemma \ref{ell on Pi_RR}] Suppose there exists some $c\in\RR_{>0}$ such that for all $g\left(\vec b\right)\in \langle \Pi_\QQ\rangle$, $$n\cdot \max_{k<n}\left|\ell\left(g\left(\vec b\right)^2\cdot\left(\left|b_{\frac{k+1}{n}}-b_{\frac{k}{n}}\right|^2\right)\right) - c\cdot \ell\left(g\left(\vec b\right)^ 2\right)\right|\longrightarrow 0\text{ as } n\rightarrow\infty.$$ Let $f_1\left(\vec b\right),\dots,f_m\left(\vec b\right)\in \langle \Pi_\QQ\rangle$ and assume, in addition, $\ell\left(g^2\left(\vec b\right)\cdot \prod_{i=1}^m{f_i}^{k_i}\left(\vec b\right)\right)\geq 0$ for all $(k_1,\dots,k_m)\in\{0,1\}^m$ and $g\in \langle \Pi_\QQ\rangle$. 
Then there exists an adapted probability space $\left(\Gamma,\left(\cG_t\right)_{t\in[0,1]},\gamma\right)$ and a process $\left(\tilde b_t\right)_{t\in[0,1]}$ with continuous paths defined on $\Gamma$ such that $\tilde b$ is a Brownian motion with respect to some measure $\eta$ on $\Gamma$ and such that not only
$f_i\left(\left(c\cdot \tilde b_t\right)_{t\in\RR\cap(0,1]}\right)\geq 0$ $\gamma$-almost surely for all $i\in\{1,\dots,m\}$, but also \begin{align*}\ell\left(p\left(\vec b\right)\right)= \EE_\gamma\left[p\left(c\cdot\vec{ \tilde b}\right)\right]
\end{align*} for all $p\left(\vec  b\right)\in\langle\Pi_\QQ\rangle$.

Moreover, there is a unique Wiener $L^1$-continuous extension from $\ell$ to $\langle\Pi_\RR\rangle$.
\end{theorem}

Similar results may, in principle, be obtained  when replacing $\QQ$ with $\QQ'=\QQ\cup F$ for a finite set of real numbers $F$. However, we shall avoid the additional notational and technical difficulties, by confining ourselves to $\QQ$.

We shall employ the theory of Loeb measures in order to apply the results for the finite-dimensional moment problem to the Wiener space and to the path space $\RR^{\QQ\cap[0,1]}$. This will enable us to prove the Theorem stated previously.

In concluding the introduction, we list a few references on (probabilistic) nonstandard analysis.

For the general theory of Loeb spaces and its pivotal extensions to applications in stochastic analysis, see Loeb's original paper \cite{momentL75} and the works by Anderson \cite{momentA76}, Hoover and Perkins \cite{momentHP83}, Keisler \cite{momentK84}, as well as Stroyan and Bayod \cite{momentSB86} -- a detailed exposition of nonstandard methods and their applications is the monograph by Albeverio {\em et al.} \cite{momentAFHL86}. 

The universality and other specific model-theoretic features of hyperfinite adapted probability spaces and their Loeb extensions when compared to ordinary filtered probability spaces, has been proven by Hoover and Keisler \cite{momentHK84} (see also Fajardo and Keisler \cite{momentFK02}).

{\em A priori}, there may be some reservations about the use of nonstandard enlargements, because the common route to construct them is via ultrapowers with respect to an ultrafilter that extends the filter of co-finite subsets of $\NN$ and whose existence comes from Zorn's Lemma. However, a construction of a definable (over {\bf ZFC}) $\aleph_1$-saturated nonstandard model of the reals has been found by Kanovei and Shelah \cite{momentS04}. Their technique has been extended by Herzberg \cite{momentH07} to allow for inductive chains of definable bounded ultrapowers, in order to prove even the existence of a definable nonstandard enlargement of the full superstructure over the reals.

\newpage
\section{Polynomials of stochastic processes}

Consider a stochastic process $(x_t)_{t\in(0,1]}$ on some filtered measurable space $\left(\Omega,\left(\cF_{t}\right)_{t\in(0,1]}\right)$. Also, let $x_0(\omega):=0$ for all $\omega\in\Omega$ (thus pinning the process $x$ to null at time zero). In introducing {\em polynomials of stochastic processes}, we simply regard a stochastic process $(x_t)_{t\in(0,1]}$ as a family of variables $x_t$, $t\in(0,1]$ --- and a polynomial of $x$ is then simply a polynomial in finitely many $x_{t_1}, \dots,x_{t_N}$ for some $N\in\NN$:

\begin{definition}\label{Pi_} Let $\Pi_\QQ$ be the following space of random variables:
\begin{align*}\Pi_\QQ&:=\left\{ Y:\Omega\rightarrow \RR \ : \ \begin{array}{c}\exists m\in\NN_0 \quad \exists i_1,\cdots,i_m\in\NN \\  \exists  q_1,\dots,q_m\in\QQ \cap(0,1]\\ Y: \omega\mapsto \prod_{j=1}^m x_{q_j}(\omega)^{i_j}\end{array}\right\}\\ &= \left\{ \prod_{j=1}^m{x_{q_j}}^{i_j} \ : \ m\in\NN_0, \quad i_1,\cdots,i_m\in\NN, \quad  q_1,\dots,q_m\in\QQ \cap(0,1]\right\}.\end{align*}

The vector space of {\em polynomials of $x$ at rational times} will be just the vector space $\langle\Pi_\QQ\rangle$ generated by $\Pi_\QQ$, whereas we define the {\em polynomials of $x$} to be the elements of the vector space $\langle\Pi_\RR\rangle$ generated by 
\begin{align*}\Pi_\RR:= \left\{ \prod_{j=1}^m{x_{t_j}}^{i_j} \ : \ m\in\NN_0, \quad i_1,\cdots,i_m\in\NN, \quad  t_1,\dots,t_m\in (0,1]\right\} .\end{align*} 

Analogously, we may define 

\begin{align*}\Pi_\aR:= &\left\{ \prod_{j=1}^m{x_{t_j}}^{i_j} \ : \ \begin{array}{c} m\in\aN_0, \quad\left( i_1,\cdots,i_m\right)\in\aN^m , \\ \left(t_1,\dots,t_m\right)\in \left(\aR\cap{{^\ast}(0,1]}\right)^m \text{ internal} \end{array} \right\} \end{align*}
(where we take $A^m$ for internal sets $A$ to be shorthand for the set of {\em internal} $m$-tuples of $A$) and
\begin{align*}\Pi_\aQ:= &\left\{ \prod_{j=1}^m{x_{t_j}}^{i_j} \ : \ \begin{array}{c} m\in\aN_0, \quad\left( i_1,\cdots,i_m\right)\in\aN^m , \\ \left(t_1,\dots,t_m\right)\in \left(\aQ\cap{^\ast}(0,1]\right)^m \text{ internal} \end{array} \right\} .\end{align*}
\end{definition}

As a another notational convention, let $\cM_1\left(\RR^d\right)$ for all $d\in\NN$ denote the space of all Borel probability measures on $\RR^d$.


\begin{lemma} \label{approxLoeb} Consider $d\in\aN$ and $\left\{t_1,\dots,t_d\right\}$ an internal subset of $\aQ\cap(0,1]$ that contains $\QQ\cap(0,1]$. Let $H\in\aN$. Define, for any $K\in\aN$ the {\em $K$-truncated $\frac{1}{H!}$-rounding operation} $\rho_{d,K,1/H!}$ to be the map $$\rho_{d,K,1/H!}:\aR^d\rightarrow\aR^d,\quad x\mapsto\left( \max\left\{y \leq x_i \ : \ y\in[-K,K]\cap\frac{1}{H!}\aZ\right\} \vee -K \right)_{i=1}^d$$ and denote by $\LL_{d,K,1/H!}$ the range of $\rho_{d,K,1/H!}$, that is the lattice $$ \LL_{d,K,1/H!}:=[-K,K]^d\cap \frac{1}{H!}\aZ^d.$$ (The subscript $d$ will be dropped where no ambiguity can arise.) Now consider any $\mu\in{^*}\cM_1\left(\RR^d\right)$, ${^*}\cM_1\left(\RR^d\right)$ being the value of the $^*$-image of the function $n\mapsto\cM_1\left(\RR^n\right)$ at $n=d\in\aN$. Suppose all $p\in\langle \Pi_\QQ\rangle$ are $^*$-integrable with respect to $\mu$. Then there exist $H,K\in\aN\setminus\NN$ such that $$\aint_{\aR^d} p\left(\vec y\right)\ \mu\left(d\vec y\right) \approx \int_{\LL_{K,1/H!}} °p \left(\ualpha \right)\ dL\left(\mu\circ\left(\rho_{K,1/H!}\right)^{-1}\right)\left(\ualpha\right)$$ 
for all $p\left(\vec x\right)\in\langle \Pi_\QQ\rangle $ for which $\aint_{\aR^d} p\left(\vec y\right)\ \mu\left(d\vec y\right)$ is finite.
\end{lemma}

\begin{proof}[Proof of Lemma \ref{approxLoeb}] It is enough to prove the Lemma for all $p\left(\vec x\right)\in\Pi_\QQ$, as it then follows for all $p\left(\vec x\right)\in\langle\Pi_\QQ\rangle$ thanks to the linearity of the integral. Recall that every finite Borel measure $\mu$ is regular (cf e.g. Bauer \cite[Lemma 26.2]{momentBauerMI}) and that any measure $\nu $ which has the positive part of a $\mu$-integrable polynomial as $\mu$-density will again be a finite Borel measure. From this we may deduce that for all $d\in\NN$, $\mu\in\cM_1\left(\RR^d\right)$, and $\Pi\subset \cL^1\left(\mu\right)$ a finite set of $\mu$-integrable polynomials, as well as for all $\eps>0$, there exists some $\bar K\in\NN$ such that for all $K\geq \bar K$ and for all $p\in\Pi$, one has $$\left|\int_{\left\{|\cdot|\leq K\right\}}p \ d\mu-\int p \ d\mu\right|<\eps $$ and such that there exists an $H\geq \bar K$ satisfying $$\forall a\in\left[-\bar K,\bar K\right)^d\cap\frac{1}{H!}\ZZ^d\quad \left| \int_{\left[a,a+\left(\frac{1}{H!}\right)_{i=1}^d\right)}p \ d\mu - p(a)\mu{\left[a,a+\left(\frac{1}{H!}\right)_{i=1}^d\right)}\right|<\eps$$ (in order to see this, first choose $\bar K$ to satisfy the first condition, and then choose $H$ depending on the choice of $\bar K$ to satisfy the second condition).
Applying the Transfer Principle to this proposition, we conclude that for our choices of $d\in\aN$ and  $\mu\in{^*}\cM_1\left(\RR^{d}\right)$, there exist for all $\varepsilon\in\aR_{>0}$, $M\in\aN$ and $\left\{p_1\left(\vec X\right),\dots,p_M\left(\vec X\right)\right\}\subset\aR\left[X_{t_1},\dots,X_{t_d}\right]$ hyperfinite numbers $H,K\in\aN\setminus\NN$ such that for $p\in\left\{p_1\left(\vec X\right),\dots,p_M\left(\vec X\right)\right\}$, one has $$\left|\aint_{\left\{|\cdot|\leq K\right\}}p \ d\mu-\aint p \ d\mu\right|<\eps$$  and $$\forall a\in\left[-K,K\right)^d\cap\frac{1}{H!}\ZZ^d\quad  \left| \aint_{\left[a,a+\left(\frac{1}{H!}\right)_{i=1}^d\right)}p \ d\mu - p(a)\mu{\left[a,a+\left(\frac{1}{H!}\right)_{i=1}^d\right)}\right|<\eps.$$
Choosing $\eps$ to be infinitesimally small, we obtain \begin{align}\aint_{\left\{|\cdot|\leq K\right\}}p \ d\mu\approx\aint p \ d\mu\label{truncate}\end{align} and \begin{align} \label{grid_integral} \forall a\in\left[-K,K\right)^d\cap\frac{1}{H!}\ZZ^d\quad  \int_{\left[a,a+\left(\frac{1}{H!}\right)_{i=1}^d\right)}p \ d\mu \approx p(a)\mu{\left[a,a+\left(\frac{1}{H!}\right)_{i=1}^d\right)} \end{align}

In particular, we may choose $M$ and $\left\{p_1,\dots,p_M\right\}$ such that $\left\{p_1\left(\vec x\right),\dots,p_M\left(\vec x\right)\right\}\supset \Pi_\QQ$, since $\Pi_\QQ$ is a subset of the (hyperfinite) set of monomials in $\aR\left[X_{\frac{1}{H!}},\dots, X_{\frac{H!}{H!}}\right]$ of total degree less than some fixed infinite hyperfinite number. 
On the other hand, thanks to the hyperfinite additivity of internal $^*$-integrals $$ \aint_{\left\{|\cdot|\leq K\right\}}p \ d\mu = \sum_{a\in\left[-K,K\right)^d\cap\frac{1}{H!}\aZ^d} \aint_{\left[a,a+\left(\frac{1}{H!}\right)_{i=1}^d\right)}p \ d\mu$$ which implies 
\begin{align} \nonumber & \left|\aint_{\left\{|\cdot|\leq K\right\}}p \ d\mu - \sum_{a\in\left[-K,K\right)^d\cap\frac{1}{H!}\aZ^d} p(a)\mu{\left[a,a+\left(\frac{1}{H!}\right)_{i=1}^d\right)}\right|\\ &\leq \mu\left\{|\cdot|\leq K\right\}\cdot\sup_{a\in\left[-K,K\right)^d\cap\frac{1}{H!}\aZ^d}\left| \aint_{\left[a,a+\left(\frac{1}{H!}\right)_{i=1}^d\right)}p \ d\mu - p(a)\mu{\left[a,a+\left(\frac{1}{H!}\right)_{i=1}^d\right)}\right|.
\label{LLdecompos}\end{align}
Note that the supremum on the right hand side of estimate \eqref{LLdecompos} is the supremum of an internal function over a hyperfinite set, and hence must be attained (i.e., it is in fact a maximum). However, the range of this internal function over this set is a set of infinitesimals for any $p$ due to condition \eqref{grid_integral}. Hence the supremum occurring on the right hand side of estimate \eqref{LLdecompos} is, indeed, the maximum of a set of infinitesimals and therefore must be infinitesimal as well. But $\mu\left\{|\cdot|\leq K\right\}\leq 1$, therefore
\begin{align}\aint_{\left\{|\cdot|\leq K\right\}}p \ d\mu \approx \sum_{a\in\left[-K,K\right)\cap\frac{1}{H!}\aZ^d} p(a)\mu{\left[a,a+\left(\frac{1}{H!}\right)_{i=1}^d\right)}.\label{intsum}
\end{align}
Observe next that $\LL_{K,1/H!}$ is hyperfinite and hence by classical Loeb integration theory, one has \begin{align*}\sum_{a\in\left[-K,K\right)\cap\frac{1}{H!}\aZ^d} p(a)\mu{\left[a,a+\left(\frac{1}{H!}\right)_{i=1}^d\right)}\\ = {\aint_{\LL_{K,1/H!}}}p\left(\ualpha\right)\ d\left(\mu\circ\left(\rho_{K,1/H!}\right)^{-1}\right)\left(\ualpha\right) \\ \approx  {\int_{\LL_{K,1/H!}}} °p\left(\ualpha\right)\ dL\left(\mu\circ\left(\rho_{K,1/H!}\right)^{-1}\right)\left(\ualpha\right) 
\end{align*}
for all $p\in\langle\Pi_\QQ\rangle$ for which the left hand side is finite (rather than merely an element of $\aR$).
Inserting this into the approximate identity \eqref{intsum} and combining it with relation \eqref{truncate} yields \begin{align}\nonumber \int_{\LL_{K,1/H!}} °p\left( \ualpha\right)\ dL\left(\mu\circ\left(\rho_{K,1/H!}\right)^{-1}\right)\left(\ualpha\right)  &\approx \aint_{\left\{|\cdot|\leq K\right\}}p \ d\mu \\ & \approx \aint_{\aR^ d}p \ d\mu.
\end{align}
Hence we have proven the identity asserted in the Lemma for all $p\left(\vec x\right)\in \Pi_\QQ$. Due to the linearity of integrals, this is sufficient to deduce the identity for all $p\left(\vec x\right)\in\langle\Pi_\QQ\rangle$.
\end{proof}

For the following remark and beyond, we shall introduce the abbreviation $$\TT_{H!}:=\left\{\frac{1}{H!},\dots, \frac{H!}{H!}\right\}.$$ 

\begin{rem} \label{PolynomialCorrespondence}We shall identify the elements of $\langle\Pi_\RR\rangle$ and $\langle \Pi_\QQ\rangle$ with elements of $\bigcup_{\vec t\in(0,1]^{<\aleph_0}}\RR\left[X_{\vec t}\right]$ and $\bigcup_{\vec t\in\left({\QQ \cap(0,1]}\right)^{<\aleph_0}}\RR\left[X_{\vec t}\right]$, respectively ($X_{\vec t}$ denoting $\left(X_{t_1},\dots,X_{t_m}\right)$ for all $m$-tuples $\vec t=\left(t_1,\dots,t_m\right)$ and any $m\in\aN_0$). Then $\ell: \langle \Pi_\QQ\rangle\rightarrow \RR$ becomes a linear map from $\bigcup_{\vec t\in\left(\QQ\cap(0,1]\right)^{<\aleph_0}}\RR\left[X_{\vec t}\right]$ (which is the vector space of polynomials in variables $X_t$ for $t\in\aQ\cap(0,1]$) to $\RR$. Therefore the ${^*}$-image ${^*}\ell$ of this map will be a map from $\bigcup\left\{\aR\left[X_{\vec t}\right] \ : \ { \vec t\in\left(\aQ\cap(0,1]\right)^{<\aleph_0} } \right\}$ to $\aR$. If we now identify elements of $\bigcup_{\vec t\in\left(\aQ\cap(0,1]\right)^{<\aleph_0} }\aR\left[X_{\vec t}\right] = \bigcup\left\{\aR\left[X_{\vec t}\right] \ : \ { \vec t\in\left(\aQ\cap(0,1]\right)^{<\aleph_0} } \right\}$ and $\langle \Pi_\aQ\rangle$, ${^*}\ell$ can be thought of as a (${^*}$-linear) map ${^*}\ell:\langle\Pi_\aQ\rangle\rightarrow\aR$.

Note for the following that the ``process'' $\vec x:=\left(x_t\right)_{t\in\TT_{H!}}$ only occurs as argument of a polynomial, and $p\left(\vec x\right)$, for $p\in\aR\left[X_{1/H!},\dots,X_{H!/H!}\right]$, will always denote an element of $\langle\Pi_\aR\rangle$.

\end{rem}
\newpage
\section{Extending Wiener $L^1$-continuous functionals from rational to irrational times}

Denote the Wiener measure on $\Omega:=C^0\left([0,1],\RR\right)\cap\left\{f\ : \ f(0)=0\right\}$ by $\WW$ and let for all $t\in[0,1]$ the random variable $b_t:\Omega\rightarrow\RR$ denote the path-space projection $b_t:\omega\mapsto \omega(t)$.

Let $\ell$ be an $L^1(\WW)$-continuous map from a superspace of the vector space $\langle \Pi_\QQ\rangle$ of {\em polynomials of Brownian motion at rational times}, which is the vector space of real-valued random variables generated by \begin{align} \nonumber \Pi_\QQ &:=\left\{ Y:C^0\left([0,1],\RR\right)\rightarrow \RR \ : \ \begin{array}{c}\exists m\in\NN_0 \quad \exists i_1,\cdots,i_m\in\NN \\ \exists  q_1,\dots,q_m\in\QQ \cap(0,1]\\ Y: \omega\mapsto \prod_{j=1}^m\omega(q_j)^{i_j}\end{array}\right\}\\ \nonumber &= \left\{ \prod_{j=1}^m{b_{q_j}}^{i_j} \ : \ m\in\NN_0, \quad i_1,\cdots,i_m\in\NN, \quad  q_1,\dots,q_m\in\QQ \cap(0,1]\right\}, \end{align}
and assume, $\ell$ has been normalised: $\ell(1)=1$.

\begin{lemma}\label{ell on Pi_RR} $\ell$ is defined on the vector space $\langle\Pi_\RR\rangle $ of {\em polynomials of Brownian motion}, which is the vector space of real-valued random variables generated by \begin{align}\Pi_\RR:= \left\{ \prod_{j=1}^m{b_{t_j}}^{i_j} \ : \ m\in\NN_0, \quad i_1,\cdots,i_m\in\NN, \quad  t_1,\dots,t_m\in (0,1]\right\} .\end{align} Moreover, if for some $f\left(\vec b\right)\in\langle\Pi_\QQ\rangle$, one has $\ell\left(g\left(\vec b\right)^2\cdot f\left(\vec b\right)\right)\geq 0$ for all $g\left(\vec b\right)\in\langle\Pi_\QQ\rangle$, then even $\ell\left(g\left(\vec b\right)^2\cdot f\left(\vec b\right)\right)\geq 0$ for all $g\left(\vec b\right)\in\langle\Pi_\RR\rangle$.
\end{lemma}
\begin{proof} Since $\Omega=C^0\left([0,1],\RR\right)\cap\left\{f\ : \ f(0)=0\right\}$ by definition, for all $t\in(0,1)$, $$\forall\omega\in\Omega\quad \lim_{s\rightarrow t} b_s(\omega) =  \lim_{s\rightarrow t} \omega(s)=\omega(t)=b_t(\omega).$$ Hence we can approximate elements of $\Pi_\RR$ by elements of $\Pi_\QQ$ pointwise on $\Omega$. Consider now any $m\in\NN$, $i_1,\cdots,i_m\in\NN$, $t_1,\dots,t_m\in (0,1]$. By applying the generalised H\"older inequality and afterwards Doob's inequality, one has a constant $C_0>0$ such that for all $\varepsilon\in(0,1)$ and $k\in\{1,\dots,m\}$ the estimates \begin{align*}& \left\| \sup_{s_k\in (t_k-\varepsilon,t_k+\varepsilon)} \left|b_{t_k}\right|^{i_k}\prod_{j\neq k}\sup_{s_j\in (t_j-\varepsilon,t_j+\varepsilon)} \left|b_{t_j}\right|^{i_j} \right\|_{L^1(\WW)} \\ &\leq \left\|\sup_{s_k\in (t_k-\varepsilon,t_k+\varepsilon)} {\left| b_{t_k}\right|}^{i_k}\right\|_{L^m(\WW)}\cdot \prod_{j\neq k}\left\|\sup_{s_j\in (t_j-\varepsilon,t_j+\varepsilon)} \left|b_{t_j}\right|^{i_j} \right\|_{L^m(\WW)}\\ &\leq C_0\cdot \sup_{s_k\in (t_k-\varepsilon,t_k+\varepsilon)} \left\|{ b_{t_k}}^{i_k}\right\|_{L^m(\WW)}\cdot \prod_{j\neq k}\left\|\sup_{s_j\in (t_j-\varepsilon,t_j+\varepsilon)} \left|b_{t_j}\right|^{i_j} \right\|_{L^m(\WW)} \end{align*} hold. This yields inductively in $k$ for all $\varepsilon\in(0,1)$ \begin{align*}& \left\| \sup_{s_j\in (t_j-\varepsilon,t_j+\varepsilon)}\prod_{j=1}^m \left|b_{s_j}\right|^{i_j} \right\|_{L^1(\WW)} \\ &\leq {C_0}^m\cdot \sup_{s_1\in (t_1-\varepsilon,t_1+\varepsilon)} \cdots\sup_{s_m\in (t_m-\varepsilon,t_m+\varepsilon)} \prod_{j=1}^m \left\|{b_{s_j}}^{i_j} \right\|_{L^m(\WW)} \\ &<+\infty\end{align*} 

Therefore if we approximate each $t_j$ by a sequence $\left(q_{n,j}\right)_{n\in\NN}$, we will be able to bound $\prod_{j=1}^m \left|b_{q_{n,j}}\right|^{i_j}$ uniformly in $n$ by a random variable in $L^1(\WW)$. But on the other hand we have already seen that $\lim_{n\rightarrow\infty}\prod_{j=1}^m \left|b_{q_{n,j}}\right|^{i_j}(\omega)=\prod_{j=1}^m \left|b_{t_{j}}(\omega)\right|^{i_j}$. Hence we get by Lebesgue's Dominated Convergence Theorem that $\lim_{n\rightarrow\infty}\prod_{j=1}^m \left|b_{q_{n,j}}\right|^{i_j}=\prod_{j=1}^m \left|b_{t_{j}}\right|^{i_j}$ in $L^1(\WW)$. Taking advantage of the assumption that $\ell$ be $L^1(\WW)$-continuously extendible, we finally conclude $$\lim_{n\rightarrow\infty}\ell\left(\prod_{j=1}^m \left|b_{q_{n,j}}\right|^{i_j}\right)=\ell\left(\prod_{j=1}^m \left|b_{t_{j}}\right|^{i_j}\right).$$ The second part of the Lemma follows from this identity.

\end{proof}

In analogy to Definition \ref{Pi_}, let us now also introduce the following notation:

\begin{definition}
\begin{align}\Pi_\aR:= &\left\{ \prod_{j=1}^m{b_{t_j}}^{i_j} \ : \ \begin{array}{c} m\in\aN_0, \quad\left( i_1,\cdots,i_m\right)\in\aN^m , \\ \left(t_1,\dots,t_m\right)\in \left(\aR\cap(0,1]\right)^m \text{ internal} \end{array} \right\} \end{align}
(where we take $A^m$ for internal sets $A$ to be shorthand for the set of {\em internal} $m$-tuples of $A$) and
\begin{align} \nonumber \Pi_\aQ:= &\left\{ \prod_{j=1}^m{b_{t_j}}^{i_j} \ : \ \begin{array}{c} m\in\aN_0, \quad\left( i_1,\cdots,i_m\right)\in\aN^m , \\  \left(t_1,\dots,t_m\right)\in \left(\aQ\cap(0,1]\right)^m \text{ internal} \end{array} \right\} .\end{align}
\end{definition}
{\em Mutandis mutatis}, viz. replacing the process $x$ by $b$, Remark \ref{PolynomialCorrespondence} about the domains of $\ell$ and $\aell$ naturally applies to this setting as well. In light of Lemma \ref{ell on Pi_RR}, through which $\ell$ is defined on all of $\langle\Pi_\RR\rangle$, and Remark \ref{PolynomialCorrespondence}, $\aell$ will be defined on all of $\Pi_\aR$ and hence, by $^*$-linearity, also on the $^*$-linear hull $\langle\Pi_\aR\rangle$ of $\Pi_\aR$.
\newpage
\section{The moment problem for linear functionals of polynomials of Brownian motion}

For all $p\left(\vec b\right)\in\Pi_\aR$ that correspond --- in the sense of Remark \ref{PolynomialCorrespondence} --- to the polynomial $p\left(\vec X\right)=\sum_k a_k\prod_{j=1}^{H!} {X_{j/H!}}^{i_j(k)}\in \aR\left[X_{\frac{1}{H!}}, \dots, X_{\frac{H!}{H!}}\right]$ and $\vec y=\left(y_{1/H!},\dots,y_{H!/H!}\right)\in\aR^{H!}$, $p\left(\vec y\right)$ will denote the hyper-real number $$ p\left(\vec y\right) = \sum_k a_k\prod_{j=1}^{H!} {y_{j/H!}}^{i_j(k)} .$$ 

In a similar vein, for all $q\left(\vec b\right)\in\Pi_\RR$ such that $q\left(\vec b\right)=\sum_k a_k\prod_{j=1}^{m(k)} {b_{t_j(k)}}^{i_j(k)}$ and $\vec y = \left(y_t\right)_{t\in\RR\cap(0,1]}$, $$ q\left(\vec y\right) = \sum_k a_k\prod_{j=1}^{m(k)} {y_{t_j(k)}}^{i_j(k)} \in\RR.$$ Also, for any $\vec y\in \aR^{{H!}}$ or $\vec y\in\RR^{\QQ\cap(0,1]}$, $y_0$ will always be zero by definition: $y_0:=0$.

Let us now come to the main results of this paper. For this sake, we shall introduce the following, fairly self-explanatory, manner of speaking:

\begin{definition} Let $c>0$ be a standard real number. We say that $\ell$ commutes with the quadratic variation scaled by $c$, if and only if $$n\cdot \max_{k<n}\left|\ell\left(g\left(\vec b\right)^2\cdot\left(\left|b_{\frac{k+1}{n}}-b_{\frac{k}{n}}\right|^2\right)\right) - c\cdot \ell\left(g\left(\vec b\right)^ 2\right)\right|\longrightarrow 0\text{ as } n\rightarrow\infty$$ and for all $g\left(\vec b\right)\in \langle \Pi_\QQ\rangle$. 
\end{definition}

\begin{lemma} \label{approxell} Suppose $\ell$ commutes with the quadratic variation scaled by some standard $c>0$. Let $f_1\left(\vec b\right),\dots,f_m\left(\vec b\right)\in \langle \Pi_\QQ\rangle$ and assume, in addition, $\ell\left(g^2\left(\vec b\right)\cdot \prod_{i=1}^m{f_i}^{k_i}\left(\vec b\right)\right)\geq 0$ for all $(k_1,\dots,k_m)\in\{0,1\}^m$ and $g\in \langle \Pi_\QQ\rangle$ (which, by Lemma \ref{ell on Pi_RR}, entails this estimate for all $g\in\langle\Pi_\RR\rangle$ as well). Then for all $c_1,c_0\in\aR$ satisfying $c_0<c<c_1$, and for all $N\in\aN$ and $H\geq N$, there exists a ${^*}$-Borel probability measure $\mu$ concentrated on $\bigcap_{k<H!}\left\{ \vec y \in \aR^{H!}\ : \ \frac{c_0}{H!}\leq \left|y_{\frac{k+1}{H!}}-y_{\frac{k}{H!}}\right|^2\leq\frac{c_1}{H!}\right\}\cap\bigcap_{i=1}^m\left\{ \vec y \in \aR^{H!}\ : \ f_i\left( \vec y \right)\geq 0\right\}$ such that not only $\ell\left(p\left(\vec b\right)\right)=\aint p\ d\mu$ for all $p\left(\vec X\right)\in \aR\left[X_{\frac{1}{H!}}, \dots, X_{\frac{H!}{H!}}\right]$, but also (by virtue of Lemma \ref{approxLoeb}) for all $p\in\langle\Pi_\QQ\rangle$, $\ell\left(p\left(\vec b\right)\right)\approx \int_{\LL_{K,1/\tilde H!}} °p\circ\rho_{K,1/\tilde H!}\ dL\left(\mu\circ\left(\rho_{K,1/\tilde H!}\right)^{-1}\right)$.
\end{lemma}
\begin{proof} Let $c_1>c>c_0>0$ in $\aR$. First, let us apply the Transfer Principle to the convergence assertion entailed by the assumption that $\ell$ commutes with quadratic version. Then there exists some $N\in\aN$ (without loss of generality, $N\in\aN\setminus \NN$) such that for all $H\geq N$, $g\in\aR\left[X_{\frac{1}{H!}}, \dots, X_{\frac{H!}{H!}}\right]$, $k<H!$, both $$\aell\left(g\left(\vec b\right)^2\cdot\left(\left|b_{\frac{k+1}{H!}}-b_{\frac{k}{H!}}\right|^2-\frac{c_0}{H!}\right)\right)\geq 0$$ and $$\aell\left(g\left(\vec b\right)^2 \cdot\left(\left|b_{\frac{k+1}{H!}}-b_{\frac{k}{H!}}\right|^2-\frac{c_1}{H!}\right)\right)\leq 0, $$ in addition to $\aell\left(g^2\left(\vec b\right)\cdot f_i\left(\vec b\right)\right)\geq 0$ for all $i\in\{1,\dots,m\}$ and for all $g\in{^*}\left(\langle \Pi_\RR\rangle\right)= \langle \Pi_\aR\rangle_{\aR}$. Now we employ the Transfer Principle in order to use Schm\"udgen's solution to the moment problem in the nonstandard universe. We intend to apply the result of this transfer to the $^*$-linear internal functional $$F:\aR\left[X_{\frac{1}{H!}}, \dots, X_{\frac{H!}{H!}}\right] \rightarrow \aR, \quad p\left(\vec X\right)\mapsto \ell\left(p\left(\vec b\right)\right),$$ as then we would find a ${^*}$-Borel measure $\mu$ with the properties asserted in the Lemma. However, due to $y_0=0$, the set $\bigcap_{k<H!}\left\{ \vec y \in \aR^{H!}\ : \ \frac{c_0}{H!}\leq \left|y_{\frac{k+1}{H!}}-y_{\frac{k}{H!}}\right|^2\leq\frac{c_1}{H!}\right\}$ is bounded by $\left(2c+\eps\right)\sqrt{H!}$ for all standard $\eps>0$ and therefore $^*$-compact. Hence, Schm\"udgen's solution to the moment problem may actually be applied in this setting, finally providing us with the $^*$-Borel measure $\mu$ as asserted in the Lemma.

\end{proof}

\begin{lemma} \label{internalWienerMoment} Let $c_1>c>c_0>0$ in $\aR$ with $c_1\approx c \approx c_0$, $c\in\RR$, and suppose $\ell$ commutes with the quadratic variation scaled by $c$. Let $f_1\left(\vec b\right),\dots,f_m\left(\vec b\right)\in \langle \Pi_\QQ\rangle$ and assume $\ell\left(g^2\left(\vec b\right)\cdot \prod_{i=1}^m{f_i}^{k_i}\left(\vec b\right)\right)\geq 0$ for all $(k_1,\dots,k_m)\in\{0,1\}^m$ and $g\in \langle \Pi_\QQ\rangle$. 

Then there is a hyperfinite $H\in\aN$, a $^*$-Borel measure $\nu$ on $\aR^{H!}\otimes 2^{H!}$, as well as an internal process $\xi:\aR^{H!}\otimes \TT_{H!}\cup\{0\}\rightarrow [c_0,c_1]$ such that 
\begin{itemize}
\item $\aell\left(p\left(\vec b\right)\right) = \aint_{\aR^{H!}\otimes\Gamma} p\left(\left(\xi_t\left(\vec\alpha\right)\cdot B_t(\omega)\right)_{t\in\TT_{H!}}\right)\ \nu\left(d\left(\vec\alpha,\omega\right)\right)  $ for all $^*$-polynomials $p\left(\vec X\right)\in \aR\left[X_{\frac{1}{H!}}, \dots, X_{\frac{H!}{H!}}\right]$, wherein
\item $B:=\left(B_t\right)_{\TT_{H!}}:\Gamma\otimes \TT_{H!}\rightarrow \frac{1}{\sqrt{H!}}\aZ$ denotes Anderson's random walk on the internal probability space $\Gamma=2^{H!}$, and
\item with $\nu$-probability $1$, $f_i\left(\left(\xi_t\cdot B_t\right)_{t\in\TT_{H!}}\right)\geq 0$ for all $i\in\{1,\dots,m\}$.
\end{itemize}
\end{lemma}
\begin{proof} According to Lemma \ref{approxell}, there is a ${^*}$-Borel probability measure $\mu$ concentrated on $\bigcap_{k<H!}\left\{\vec y\in \aR^{H!}\ : \ \frac{c_0}{H!}\leq \left|y_{\frac{k+1}{H!}}-y_{\frac{k}{H!}}\right|^2\leq\frac{c_1}{H!}\right\}\cap\bigcap_{i=1}^m\left\{\vec y\in \aR^{H!}\ : \ f_i\left(\vec y\right)\geq 0\right\}$ such that $\ell\left(p\left(\vec b\right)\right)=\aint p\ d\mu$ for all $p\left(\vec X\right)\in \aR\left[X_{\frac{1}{H!}}, \dots, X_{\frac{H!}{H!}}\right]$. Now we construct a transformation $\psi:\vec w \mapsto\left( \xi_t\left(\vec w\right), B_t\left(\varphi\left(\vec w\right)\right)\right)_{t\in\TT_{H!}}$ for some internal process $\xi:\aR^{H!}\otimes\left(\TT_{H!}\cup\{0\}\right)\rightarrow [c_0,c_1]$, a map $\varphi:\aR\rightarrow 2^{H!}$ and a symmetric ${^*}$-random walk $B:2^{H!}\otimes\TT_{H!}\rightarrow \frac{1}{\sqrt{H!}}\aZ$ such that $$ \forall t\in\TT_{H!}\quad \vec w =\xi_t\left(\vec w\right)\cdot B_t\left(\varphi\left(\vec w\right)\right).$$
For this sake, we simply define $\xi\left(\vec w,t\right)$ for all $\vec w\in \bigcap_{k<H!}\left\{ \vec y \in \aR^{H!}\ : \ \frac{c_0}{H!}\leq \left|y_{\frac{k+1}{H!}}-y_{\frac{k}{H!}}\right|^2\leq\frac{c_1}{H!}\right\}$ and $t\in-\frac{1}{H!}+\TT_{H!}$ b$$\xi\left(\vec w,t\right)= \frac{ \left|w_{t+\frac{1}{H!}}-w_{t}\right|}{\sqrt{H!}}$$ as well as $$\bar B\left(\vec w, t\right) := \frac{1}{\sqrt{H!}}\cdot\frac{w_{t+\frac{1}{H!}}-w_{t}}{\xi\left(\vec w,t\right)}.$$ Then for all $\vec w\in\bigcap_{k<H!}\left\{ \vec y \in \aR^{H!}\ : \ \frac{c_0}{H!}\leq \left|y_{\frac{k+1}{H!}}-y_{\frac{k}{H!}}\right|^2\leq\frac{c_1}{H!}\right\} $ one has $$ \forall t\in -\frac{1}{H!}+\TT_{H!}\quad\bar B\left(\vec w,t+\frac{1}{H!}\right)-\bar B\left(\vec w,t\right)\in\left\{\pm\frac{1}{\sqrt{H!}}\right\},$$ and therefore, there is a surjective map $\varphi:\aR^{H!}\rightarrow 2^{H!}$ and a process $B:2^{H!}\otimes \TT_{H!}\cup\{0\}\rightarrow \frac{1}{\sqrt{H!}}\aZ$ such that for all $\vec w$, $$ \forall t\in \TT_{H!}\cup\{0\}\quad \bar B\left(\vec w,t\right) = B\left(\varphi\left(\vec w\right),t\right),$$ and $B$ is just Anderson's random walk (denoted $\chi$ in Anderson's original paper \cite{momentA76}) associated to the mesh size $H!$.  Since by hyperfinite induction in $\TT_{H!}$, one can prove $$ \forall t\in \TT_{H!}\quad w_t= \xi\left(\vec w,t\right)\cdot \bar B\left(\vec w,t\right),$$ we finally obtain $$ \forall t\in \TT_{H!}\cup\{0\}\quad w_t= \xi\left(\vec w,t\right)\cdot \underbrace{\bar B\left(\vec w,t\right)}_{= B\left(\varphi\left(\vec w\right),t\right)}$$ for all $\vec w\in\bigcap_{k<H!}\left\{ \vec y \in \aR^{H!}\ : \ \frac{c_0}{H!}\leq \left|y_{\frac{k+1}{H!}}-y_{\frac{k}{H!}}\right|^2\leq\frac{c_1}{H!}\right\} $. Hence we have constructed an injective transformation \begin{align*}\psi: \bigcap_{k<H!}\left\{ \vec y \in \aR^{H!}\ : \ \frac{c_0}{H!}\leq \left|y_{\frac{k+1}{H!}}-y_{\frac{k}{H!}}\right|^2\leq\frac{c_1}{H!}\right\} & \rightarrow \aR^{H!}\otimes 2^{H!},\\ \vec w & \mapsto\left( \xi_t\left(\vec w\right), B_t\left(\varphi\left(\vec w\right)\right)\right)_{t\in\TT_{H!}},\end{align*} and the image of the $^*$-Borel probability measure $\mu$ under this transformation will again be a $^*$-Borel probability measure -- this time on $\aR^{H!}\otimes 2^{H!}\subseteq \aR^{H!}\otimes \aR^{H!}=\aR^{2H!}$ and concentrated on $\psi\left[\bigcap_{k<H!}\left\{ \vec y \in \aR^{H!}\ : \ \frac{c_0}{H!}\leq \left|y_{\frac{k+1}{H!}}-y_{\frac{k}{H!}}\right|^2\leq\frac{c_1}{H!}\right\}\right]$.

\end{proof}

\begin{corollary} \label{externalWiener} Suppose again $\ell$ commutes with the quadratic variation scaled by some standard $c>0$. Let $f_1\left(\vec b\right),\dots,f_m\left(\vec b\right)\in \langle \Pi_\QQ\rangle$ and assume, in addition, $\ell\left(g^2\left(\vec b\right)\cdot \prod_{i=1}^m{f_i}^{k_i}\left(\vec b\right)\right)\geq 0$ for all $(k_1,\dots,k_m)\in\{0,1\}^m$ and $g\in \langle \Pi_\QQ\rangle$. 
Then there exists an adapted probability space $\left(\Gamma,\left(\cG_t\right)_{t\in[0,1]},\gamma\right)$ and a process $\left(\tilde b_t\right)_{t\in[0,1]}$ with continuous paths defined on $\Gamma$ such that $\tilde b$ is a Brownian motion with respect to some measure $\eta$ on $\Gamma$ and such that not only
$f_i\left(\left(c\cdot \tilde b_t\right)_{t\in\RR\cap(0,1]}\right)\geq 0$ $\gamma$-almost surely for all $i\in\{1,\dots,m\}$, but also \begin{align*}\ell\left(p\left(\vec b\right)\right)= \EE_\gamma\left[p\left(c\cdot\vec \tilde b\right)\right]
\end{align*} for all $p\left(\vec b\right)\in\langle\Pi_\QQ\rangle$.
\end{corollary}
\begin{proof} The internal measure $\nu$ constructed in Lemma \ref{internalWienerMoment} is a $^*$-Borel measure (on $\aR^{2H!}$) and therefore we may approximate the $^*$-integral of every function $p\in\langle\Pi_\QQ\rangle$ with respect to the internal measure $\nu$ by the standard integral of $p$ with respect to the Loeb extension of the hyperfinite internal measure $\nu\circ{\rho_{ 2H!,K,1/\tilde H! }}^{-1}$ for some $K\in\aN\setminus\NN$ uniformly in $p$. (This was proven in Lemma \ref{approxLoeb}.) Therefore, using the specific choice of $\nu$ and $\xi$ according to Lemma \ref{internalWienerMoment}, \begin{align}\nonumber \ell\left(p\left(\vec b\right)\right)& = °\left(\aell\left(p\left(\vec b\right)\right)\right)\\ & =  °\aint_{\aR^{H!}\otimes\Gamma} p\left(\left(\xi_t\left(\vec\alpha\right)\cdot B_t(\omega)\right)_{t\in\TT_{H!}}\right)\ \nu\left(d\left(\vec\alpha,\omega\right)\right) \\ &=  \int_{\LL_{2H!,K,1/\tilde H!}\otimes 2^{H!}} °p\left(\left(\xi_{t}\left(\vec \alpha\right)\cdot B_t(\omega)\right)_{t\in\TT_{H!}}\right)\ dL\left( \nu\circ{\rho_{2 H!,K,1/\tilde H! }}^{-1}\right)\left(\vec \alpha,\omega\right) \label{approx_ell_internally}\end{align} for all $p\in \langle \Pi_\QQ\rangle$. Note, however, that $p\left(\vec y\right)$, for any $\vec y\in \aR^{\TT_{H!}}$ only depends on the rational coordinates. 
If $L\left( \nu\circ{\rho_{2 H!,K,1/\tilde H! }}^{-1}\right)\left\{°B_t\not\in\RR\right\}>0$, then for quadratic $p\left(\vec X\right)$, identity \eqref{approx_ell_internally} would yield $\ell\left(p\left(\vec b\right)\right)=+\infty$, contradicting $\ell\left(p\left(\vec b\right)\right)\in\RR$ for all $p\in\langle\Pi_\QQ\rangle$. Therefore the event $$ \Omega':=\bigcap_{t\in\QQ\cap[0,1]}\left\{°B_t\in\RR\right\}\subseteq \LL_{2H!,K,1/\tilde H!}\otimes 2^{H!}$$ is almost certain: $L\left( \nu\circ{\rho_{2 H!,K,1/\tilde H! }}^{-1}\right)\left[\Omega'\right]=1$. Also, we recall that the standard part $\tilde b$ of $B$ has been defined by Anderson \cite[Notation 25]{momentA76} in such a way that for all (rather than merely for Wiener/Anderson-almost all) elements $\omega$ of the projection of $\Omega'$ to  $2^{H!}$ and for arbitrary $t\in\QQ\cap(0,1]$, the identity $\tilde b_t(\omega)=°B_t(\omega)$ holds, and as Anderson \cite[Theorem 26]{momentA76} has shown, $\left(\tilde b_t\right)_{t\in\RR\cap[0,1]}$ is a normalised Brownian motion.
This yields, thanks to the $S$-continuity of $p\in\langle\Pi_\QQ\rangle$, \begin{align*}\ell\left(p\left(\vec b\right)\right) = & \int_{\LL_{2H!,K,1/\tilde H!}\otimes 2^{H!}} p\left(\left(c\cdot\tilde b_t(\omega)\right)_{t\in\RR\cap(0,1]}\right)\ dL\left( \nu\circ{\rho_{ 2H!,K,1/\tilde H! }}^{-1}\right)\left(\vec \alpha,\omega\right). \end{align*} 
Now, the integrand on the right hand side does not depend on $\vec \alpha$ anymore. Hence, the projection of $L\left( \nu\circ{\rho_{ 2H!,K,1/\tilde H! }}^{-1}\right)$ to the second component $2^{H!}=\left(\LL_{2H!,K,1/\tilde H!}\otimes 2^{H!}\right)_2$ may serve as our measure $\gamma$ on $\Gamma=2^{H!}$.
\end{proof}


\newpage
\begin{thebibliography}{20}


\bibitem{momentAFHL86}
S.~Albeverio, R.~H{\o}egh-Krohn, J.E. Fenstad, and T.~Lindstr{\o}m.
\newblock {\em {Nonstandard methods in stochastic analysis and mathematical
  physics.}}
\newblock {Pure and Applied Mathematics. 122. Orlando etc.: Academic Press},
  1986.

\bibitem{momentA76}
R.M. Anderson.
\newblock {A non-standard representation for Brownian motion and Ito
  integration.}
\newblock {\em Isr. J. Math.}, 25:15--46, 1976.

\bibitem{momentBauerMI}
H.~Bauer.
\newblock {\em {Ma{\ss}- und Integrationstheorie. 2., \"uberarb. Aufl.}}
\newblock {de Gruyter Lehrbuch. Berlin: Walter de Gruyter}, 1992.


\bibitem{momentBW01}
R.~Berr and T.~W{\"o}rmann.
\newblock {Positive polynomials and tame preorderings.}
\newblock {\em Math. Z.}, 236(4):813--840, 2001.


\bibitem{momentB68}
P.~Billingsley.
\newblock {\em {Convergence of probability measures. 2nd ed.}}
\newblock {Wiley Series in Probability and Statistics. Chichester: Wiley},
  1999.
  

\bibitem{momentCh69}
G.~Choquet.
\newblock {\em {Lectures on analysis. I: Integration and topological vector
  spaces. Edited by J. Marsden, T. Lance and S. Gelbart. 3rd printing.}}
\newblock {Mathematics Lecture Note Series. 24. Reading, MA: W.A. Benjamin},
  1976.


\bibitem{momentFK02}
S.~Fajardo and H.J. Keisler.
\newblock {\em {Model theory of stochastic processes.}}
\newblock {Lecture Notes in Logic. 14. Urbana, IL: Association for Symbolic
  Logic. Natick, MA: A K Peters}, 2002.


\bibitem{momentH36}
E.K. Haviland.
\newblock {On the momentum problem for distribution functions in more than one
  dimension. II.}
\newblock {\em Am. J. Math.}, 58:164--168, 1936.

\bibitem{momentH07}
F.S. Herzberg.
\newblock {A definable nonstandard enlargement}.
\newblock Submitted.


\bibitem{momentHK84}
D.N. Hoover and H.J. Keisler.
\newblock {Adapted probability distributions}.
\newblock {\em Trans. Am. Math. Soc.}, 286:159--201, 1984.

\bibitem{momentHP83}
D.N. Hoover and E.~Perkins.
\newblock {Nonstandard construction of the stochastic integral and applications
  to stochastic differential equations. I., II.}
\newblock {\em Trans. Am. Math. Soc.}, 275:1--58, 1983.

\bibitem{momentS04}
V.~Kanovei and S.~Shelah.
\newblock {A definable nonstandard model of the reals.}
\newblock {\em J. Symb. Log.}, 69(1):159--164, 2004.

\bibitem{momentK84}
H.J. Keisler.
\newblock {An infinitesimal approach to stochastic analysis.}
\newblock {\em Mem. Am. Math. Soc.}, 297, 1984.

\bibitem{momentKR05}
O.~Kounchev and H.~Render.
\newblock {Pseudopositive multivariate moment problem.}
\newblock {\em C. R. Acad. Bulg. Sci.}, 58(11):1243--1246, 2005.


\bibitem{momentKM}
S.~Kuhlmann and M.~Marshall.
\newblock {Positivity, sums of squares and the multi-dimensional moment
  problem.}
\newblock {\em Trans. Am. Math. Soc.}, 354(11):4285--4301, 2002.

\bibitem{momentKMSch}
S.~Kuhlmann, M.~Marshall, and N.~Schwartz.
\newblock {Positivity, sums of squares and the multi-dimensional moment
  problem. II.}
\newblock {\em Adv. Geom.}, 5(4):583--606, 2005.


\bibitem{momentL75}
P.A. Loeb.
\newblock {Conversion from nonstandard to standard measure spaces and
  applications in probability theory.}
\newblock {\em Trans. Am. Math. Soc.}, 211:113--122, 1975.


\bibitem{momentM00}
M.~Marshall.
\newblock {\em {Positive polynomials and sums of squares}}.
\newblock {Pisa: Istituti Editoriali e Poligrafici Internazionali}, 2000.



\bibitem{momentSchm91}
K.~Schm{\"u}dgen.
\newblock {The $K$-moment problem for compact semi-algebraic sets.}
\newblock {\em Math. Ann.}, 289(2):203--206, 1991.

\bibitem{momentSB86}
K.D. Stroyan and J.M. Bayod.
\newblock {\em {Foundations of infinitesimal stochastic analysis.}}
\newblock {Studies in Logic and the Foundations of Mathematics. 119.
  Amsterdam-New York-Oxford: North-Holland}, 1986.

\bibitem{momentW98}
T.~W{\"o}rmann.
\newblock {\em Strikt positive Polynome in der semialgebraischen Geometrie}.
\newblock Dissertation, Universit{\"at} Dortmund, 1998.

\end{thebibliography}
\end{document}